\documentclass[10pt,twoside]{article}
\usepackage{Latex-document}
\usepackage{latexsym}
\newtheorem{DE}{Definition}[section]
\newtheorem{LE}[DE]{Lemma}

\newtheorem{CO}[DE]{Corollary}
\newtheorem{THTH}[DE]{Theorem}
\newtheorem{PR}[DE]{Proposition}
\newtheorem{CN}[DE]{Conjecture}

\newcommand {\sm} {\setminus}

\renewcommand{\thesection}{\arabic{section}}
\renewcommand{\thesubsection}{\thesection.\arabic{subsection}.}
\def\addsec{\addtocounter{section}{1} \setcounter{DE}{0}}

\title{\bf The Strong Perfect Graph Conjecture\thanks{This work was supported in
part by NSF grant DMI-0098427 and ONR grant N00014-97-1-0196.}\vskip 6mm}
\author{G\'{e}rard Cornu\'{e}jols\vspace*{-0.5cm}\thanks{GSIA, Carnegie Mellon University,
Schenley Park, Pittsburgh, PA 15213, USA. E-mail: gc0v@andrew.cmu.edu}}
\date{\vspace{-8mm}}

\begin{document}

\maketitle

\thispagestyle{first} \setcounter{page}{547}

\begin{abstract}

\vskip 3mm

A graph is {\em perfect} if, in all its induced subgraphs, the size of a largest clique is equal to the chromatic
number. Examples of perfect graphs include bipartite graphs, line graphs of bipartite graphs and the complements
of such graphs. These four classes of perfect graphs will be called {\em basic}. In 1960, Berge formulated two
conjectures about perfect graphs, one stronger than the other. The weak perfect graph conjecture, which states
that a graph is perfect if and only if its complement is perfect, was proved in 1972 by Lov\'asz. This result is
now known as the perfect graph theorem. The strong perfect graph conjecture (SPGC) states that a graph is perfect
if and only if it does not contain an odd hole or its complement. The SPGC has attracted a lot of attention. It
was proved recently (May 2002) in a remarkable sequence of results by Chudnovsky, Robertson, Seymour and Thomas.
The proof is difficult and, as of this writing, they are still checking the details. Here we give a flavor of the
proof. Let us call {\em Berge graph} a graph that does not contain an odd hole or its complement. Conforti,
Cornu\'ejols, Robertson, Seymour, Thomas and Vu\v{s}kovi\'c (2001) conjectured a structural property of Berge
graphs that implies the SPGC: Every Berge graph $G$ is basic or has a skew partition or a homogeneous pair, or $G$
or its complement has a 2-join.
 A {\em skew partition} is a partition of the vertices
into nonempty sets $A, B, C, D$ such that every vertex of $A$ is adjacent to every vertex of $B$ and there is no
edge between $C$ and $D$.
 Chv\'atal introduced this concept in 1985 and
conjectured that no minimally imperfect graph has a skew partition. This conjecture was proved recently by
Chudnovsky and Seymour (May 2002). Cornu\'ejols and Cunningham introduced 2-joins in 1985 and showed that they
cannot occur in a minimally imperfect graph different from an odd hole. Homogeneous pairs were introduced in 1987
by Chv\'atal and Sbihi, who proved that they cannot occur in minimally imperfect graphs. Since skew partitions,
2-joins and homogeneous pairs cannot occur in minimally imperfect Berge graphs, the structural property of Berge
graphs stated above implies the SPGC. This structural property was proved: (i) When $G$ contains the line graph of
a bipartite subdivision of a 3-connected graph (Chudnovsky, Robertson, Seymour and Thomas (September 2001)); (ii)
When $G$ contains a stretcher (Chudnovsky and Seymour (January 2002)); (iii) When $G$ contains no proper wheels,
stretchers or their complements (Conforti, Cornu\'ejols and Zambelli (May 2002));  (iv) When $G$ contains a proper
wheel, but no stretchers or their complements (Chudnovsky and Seymour (May 2002)). (ii), (iii) and (iv) prove the
SPGC.

\vskip 4.5mm

\end{abstract}

\vskip 12mm

\section*{1. Introduction} \addsec

\vskip-5mm \hspace{5mm}

In this paper, all graphs are simple (no loops or multiple edges) and finite. The vertex set of graph $G$ is
denoted by $V(G)$ and its edge set by $E(G)$. A {\em stable set} is a set of vertices no two of which are
adjacent. A {\em clique} is a set of vertices every pair of which are adjacent. The cardinality of a largest
clique in graph $G$ is denoted by $\omega (G)$. The cardinality of a largest stable set is denoted by $\alpha
(G)$. A {\em $k$-coloring} is a partition of the vertices into $k$ stable sets (these stable sets are called {\em
color classes}). The {\em chromatic number} $\chi (G)$ is the smallest value of $k$ for which there exists a
$k$-coloring. Obviously,  $\omega (G) \leq \chi (G)$ since the vertices of a clique must be in distinct color
classes of the $k$-coloring. An {\em induced subgraph} of $G$ is a graph with vertex set $S \subseteq V(G)$ and
edge set comprising all the edges of $G$ with both ends in $S$. It is denoted by $G(S)$. The graph $G(V(G) - S)$
is denoted by $G \setminus S$. A graph $G$ is {\em perfect} if
 $\omega (H) = \chi (H)$ for every induced subgraphs $H$
of $G$. A graph is {\em minimally imperfect} if it is not perfect but all its  proper induced subgraphs are.

A {\em hole} is a  graph induced by a chordless cycle of length at least $4$. A hole is {\em odd} if it contains
an odd number of vertices. Odd holes are not perfect since their chromatic number is 3 whereas the size of their
largest clique is 2. It is easy to check that odd holes are minimally imperfect. The complement of a graph $G$ is
the graph $\bar G$ with the same vertex set as $G$, and $uv$ is an edge of $\bar G$ if and only if it is not an
edge of $G$. The odd holes and their complements are the only known minimally imperfect graphs. In 1960 Berge
\cite{bp} proposed the following conjecture, known as the {\it Strong Perfect Graph Conjecture}.

\begin{CN} \label{spgc} {\bf (Strong Perfect Graph Conjecture)}
{\em (Berge \cite{bp})} The only minimally imperfect graphs are the odd holes and their complements. \end{CN}

At the same time, Berge also made a weaker conjecture, which states that a graph $G$ is perfect if and only if its
complement $\bar G$ is perfect. This conjecture was proved by Lov\'asz \cite{lo2} in 1972 and is known as the {\em
Perfect Graph Theorem}.

\begin{THTH} {\bf (Perfect Graph Theorem)}
{\em (Lov\'asz \cite{lo2})} Graph $G$ is perfect if and only if graph ${\bar G}$ is perfect. \label{L}
\end{THTH}
\noindent {\bf Proof:} Lov\'asz \cite{lo3} proved the following
stronger result.
\\
{\bf Claim 1:} {\it A graph $G$ is perfect if and only if, for
every induced subgraph $H$, the number of vertices of $H$ is at
most $\alpha (H) \omega (H)$.}

Since $\alpha (H) =  \omega (\bar{H})$ and $\omega (H) =  \alpha (\bar{H})$, Claim~1 implies Theorem~\ref{L}.
\\
{\bf Proof of Claim 1:} We give a proof of this result due to
Gasparyan \cite{gas}. First assume that $G$ is perfect. Then, for
every induced subgraph $H$, $\omega (H) = \chi (H)$. Since the
number of vertices of $H$ is at most $\alpha (H) \chi (H)$, the
inequality follows.

Conversely, assume that $G$ is not perfect. Let $H$ be a minimally imperfect subgraph of $G$ and let $n$ be the
number of vertices of $H$. Let $\alpha=\alpha(H)$ and $\omega=\omega(H)$. Then $H$ satisfies
$$ \omega =\chi (H\backslash v) \mbox{ for every vertex } v\in V(H)$$
$$\mbox{and }\;\omega = \omega (H\backslash S) \mbox{ for every
 stable set } S\subseteq V(H).$$
Let $A_0$ be an $\alpha$-stable set of $H$. Fix an $\omega $-coloring of each of the $\alpha $ graphs $H\backslash
s$ for $s\in A_0$, let $A_1,\ldots ,A_{\alpha \omega}$ be the stable sets occuring as a color class in one of
these colorings and let  ${\cal A}:=\{A_0,A_1,\ldots ,A_{\alpha \omega}\}$. Let ${\bf A}$ be the corresponding
stable set versus vertex incidence matrix. Define ${\cal B}:=\{B_0,B_1,\ldots ,B_{\alpha \omega}\}$ where $B_i$ is
an $\omega$-clique of $H\backslash A_i$. Let ${\bf B}$ be the corresponding clique versus vertex incidence matrix.
\\

{\bf Claim 2:} {\it Every $\omega$-clique of $H$ intersects all
but one of the stables sets in ${\cal A}$.}
\\

{\bf Proof of Claim 2:} Let $S_1,\ldots,S_{\omega}$ be any
$\omega$-coloring of $H \setminus v$. Since any $\omega$-clique
$C$ of $H$ has at most one vertex in each $S_i$, $C$ intersects
all $S_i$'s if $v \not\in C$ and all but one if $v \in C$. Since
$C$ has at most one vertex in $A_0$, Claim~2 follows.

In particular, it follows that ${\bf A}{\bf B}^T=J-I$. Since $J-I$ is nonsingular, ${\bf A}$ and ${\bf B}$ have at
least as many columns as rows, that is $n \geq \alpha \omega +1$. This completes the proof of Claim~1.

\section*{2. Four Basic Classes of Perfect Graphs} \addsec

\vskip-5mm \hspace{5mm}

Bipartite graphs are perfect since, for any induced subgraph $H$, the bipartition implies that $\chi (H) \leq 2$
and therefore $\omega (H) = \chi (H)$.

A graph $L$ is the {\em line graph} of a graph $G$ if $V(L)=E(G)$ and two vertices of $L$ are adjacent if and only
if the corresponding edges of $G$ are adjacent.

\begin{PR}
Line graphs of bipartite graphs are perfect.
\end{PR}
\noindent {\bf Proof:} If $G$ is bipartite, $\chi '(G)=\Delta
(G)$ by a theorem of K\"onig \cite{ko}, where $\chi '$ denotes the
edge-chromatic number and $\Delta$ the largest vertex degree.

If $L$ is the line graph of a bipartite graph $G$, then $\chi (L) = \chi '(G)$ and $\omega (L) = \Delta (G)$.
Therefore $\chi (L) = \omega (L)$. Since induced subgraphs of $L$ are also line graphs of bipartite graphs, the
result follows.

Since bipartite graphs and line graphs of bipartite graphs are perfect, it follows from Lov\'asz's perfect graph
theorem (Theorem~\ref{L}) that the complements of bipartite graphs and of line graphs of bipartite graphs are
perfect. This can also be verified directly, without using the perfect graph theorem. To summarize, in this
section we have introduced four basic classes of perfect graphs:
\begin{itemize}
\item bipartite graphs and their complements, and
\item line graphs of bipartite graphs and their complements.
\end{itemize}

\section*{3. 2-Join} \addsec

\vskip-5mm \hspace{5mm}

A graph $G$ has a {\em 2-join} if its vertices can be partitioned into sets $V_1$ and $V_2$, each of cardinality
at least three, with nonempty disjoint subsets $A_1,B_1 \subseteq V_1$ and $A_2,B_2 \subseteq V_2$, such that all
the vertices of $A_1$ are adjacent to all the vertices of $A_2$, all the vertices of $B_1$ are adjacent to all the
vertices of $B_2$ and these are the only adjacencies between $V_1$ and $V_2$. There is an $O(|V(G)|^2|E(G)|^2)$
algorithm to find whether a graph $G$ has a 2-join \cite{cc}.

When $G$ contains a 2-join, we can decompose $G$ into two blocks $G_1$ and $G_2$ defined as follows.

\begin{DE}
If $A_2$ and $B_2$ are in different connected components of $G(V_2)$, define {\em block} $G_1$ to be $G(V_1 \cup
\{ p_1, q_1 \})$, where $p_1 \in A_2$ and $q_1 \in B_2$. Otherwise, let $P_1$ be a shortest path from $A_2$ to
$B_2$ and define {\em block} $G_1$ to be $G(V_1 \cup V(P_1))$. Block $G_2$ is defined similarly.
\end{DE}

 Next we show that the 2-join decomposition
preserves perfection (Cornu\'ejols and Cunningham \cite{cc}; see also Kapoor \cite{ak} Chapter 8). Earlier, Bixby
\cite{bix} had shown that the simpler join decomposition preserves perfection.

\begin{THTH} \label{2j}
Graph $G$ is perfect if and only if its blocks $G_1$ and $G_2$ are perfect.
\end{THTH}

\noindent {\bf Proof:} By definition, $G_1$ and $G_2$ are induced
subgraphs of $G$. It follows that, if $G$ is perfect, so are
$G_1$ and $G_2$. Now we prove the converse: If $G_1$ and $G_2$
are perfect, then so is $G$. Let $G^*$ be an induced subgraph of
$G$. We must show
$$(*)  \;\;\;\;\; \omega (G^*) = \chi (G^*).$$

For $i=1,2$, let $V_i^* = V_i \cap V(G^*)$. The proof of $(*)$ is based on a coloring argument, combining $\omega
(G^*)$-colorings of the perfect graphs $G(V_1^*)$ and $G(V_2^*)$ (Claim~3) into an $\omega (G^*)$-coloring of
$G^*$ (Claim 4). To prove Claim~3, we will use the following results.

\noindent {\bf Claim~1:} { (Lov\'asz's Replication Lemma
\cite{lo2})} {\it Let $\Gamma$ be a perfect graph and $v \in
V(\Gamma)$. Create a new vertex $v'$ adjacent to $v$ and to all
the neighbors of $v$. Then the resulting graph $\Gamma '$ is
perfect.}
\\
\noindent {\bf Proof of Claim~1:} It suffices to show that $\omega
(\Gamma ') = \chi (\Gamma ')$ since, for induced subgraphs, the
proof follows similarly. We distinguish two cases. Suppose first
that $v$ is contained in some $\omega (\Gamma)$-clique of
$\Gamma$. Then $\omega (\Gamma ') = \omega (\Gamma) +1$. Since at
most one new color is needed in $\Gamma '$, $\omega (\Gamma ') =
\chi (\Gamma ')$ follows.

Now suppose that $v$ is not contained in any $\omega (\Gamma)$-clique of $\Gamma$. Consider any $\omega
(\Gamma)$-coloring of $\Gamma$ and let $A$ be the color class containing $v$. Then, $\omega (\Gamma \setminus (A -
\{ v \} )) = \omega (\Gamma) -1$, since every $\omega (\Gamma)$-clique of $\Gamma$ meets $A - \{ v \}$. By the
perfection of $\Gamma$, the graph $\Gamma \setminus (A - \{ v \} )$ can be colored with $\omega (\Gamma) - 1$
colors. Using one additional color for the vertices $(A - \{ v \} ) \cup \{ v' \}$, we obtain an $\omega
(\Gamma)$-coloring of $\Gamma '$. This proves Claim~1.

We say that $\Gamma '$ is obtained from $\Gamma$ by {\em replicating} $v$. Replication can be applied recursively.
We say that $v$ is {\em replicated $k$ times} if $k$ copies of $v$ are made, including $v$.
\\

\noindent {\bf Claim~2:}  {\it Let $\Gamma$ be a graph and $uv$ an
edge of $\Gamma$ such that the vertices $u$ and $v$ have no common
neighbor. Let $\Gamma '$ be the graph obtained from $\Gamma$ by
replicating vertex $v$ into $v'$. Let $H$ be the graph obtained
from $\Gamma '$ by deleting edge $uv'$. Then $\Gamma$ is perfect
if and only if $H$ is perfect.}
\\

\noindent {\bf Proof of Claim~2:}   If $H$ is perfect, then so is
$\Gamma$ since $\Gamma$ is an induced subgraph of $H$.

Conversely, suppose that $\Gamma$ is perfect and $H$ is not. Let $H^*$ be a minimally imperfect subgraph of $H$.
Let $\Gamma ^*$ be the subgraph of $\Gamma '$ induced by the vertices of $H^*$. Since $\Gamma ^*$ is perfect but
$H^*$ is not, $V(H^*)$ must contain vertices $u$ and $v'$. Also $\chi (\Gamma ^*) = \chi (H^*)$ and $\omega
(\Gamma ^*) = \omega (H^*)+1$. Therefore $uvv'$ is the unique maximum clique in $\Gamma ^*$ and $\omega (H^*) =
2$. The only neighbor of $v$ in $H^*$ is $u$ since otherwise $v, v'$ would be in a clique of cardinality three in
$H^*$.
 Now $v' $ is a vertex of degree 1 in $H^*$,
a contradiction to the assumption that $H^*$ is minimally imperfect. This proves Claim~2.

For $i=1,2$, let $A_i^*=A_i \cap V(G^*)$, $B_i^*=B_i \cap V(G^*)$, $a_i = \omega (A_i^*)$ and $b_i = \omega
(B_i^*)$. Let $G_i^* = G_i \setminus (V_i-V_i^*)$ and $\omega \geq \omega (G_i^*)$. In an $\omega$-coloring of
$G_i^*$, let $C(A_i^*)$ and $C(B_i^*)$ denote the sets of colors in $A_i^*$ and $B_i^*$ respectively.

\noindent {\bf Claim 3:} {\it There exists an $\omega$-coloring of
$V_i^*$ such that $|C(A_i^*)| = a_i$ and  $|C(B_i^*)| = b_i$.
Furthermore, if $G_i$ contains path $P_i$ and

(i) if $P_i$ has an odd number of edges, then $|C(A_i^*) \cap C(B_i^*)| = \max (0, a_i + b_i - \omega )$,

(ii) if $P_i$ has an even number of edges, then $|C(A_i^*) \cap
C(B_i^*)| = \min (a_i, b_i )$.}
\\

\noindent{\bf Proof of Claim~3:} First assume that block $G_i$ is
induced by $V_i \cup \{ p_i, q_i \}$. In $G_i^*$, replicate $p_i$
$\omega - a_i$ times and $q_i$ $\omega - b_i$ times. By Claim~1,
this new graph $H$ is perfect and $\omega (H) = \omega$.
Therefore an $\omega$-coloring of $H$ exists. This coloring
induces an $\omega$-coloring of $V_i^*$ with $|C(A_i^*)|=a_i$ and
$|C(B_i^*)|=b_i$. Now assume that $G_i$ contains path $P_i$. We
consider two cases.

(i) $P_i$ has an odd number of edges.

Let $P_i = x_1, \ldots, x_{2k}$. In $G_i^*$, replicate vertex $x_{2k}$ into $x'_{2k}$ and remove edge
$x_{2k-1}x'_{2k}$. By Claim~2, the new graph is perfect. For $i$ odd, $1 \leq i < 2k$, replicate vertex $x_i$
$\omega - a_i$ times. For $i$ even, $1 < i \leq 2k-2$, replicate vertex $x_i$ $a_i$ times.

If $a_i + b_i < \omega$, replicate $x_{2k}$ $a_i$ times and replicate $x'_{2k}$ $\omega - a_i - b_i$ times. By
Claim~1, this new graph $H$ is perfect. Since $\omega (H) = \omega$, $H$ has an $\omega$-coloring. Note that
$|C(A_i^*)|=a_i$ and $|C(B_i^*)|=b_i$ and every vertex of $P_i$ belongs to two cliques of size $\omega$. So the
colors that appear in the replicates of $x_{2k}$ are precisely $C(A_i)$. Therefore $B_i^*$ is colored with colors
that do not appear in $C(A_i^*)$. Thus $|C(A_i^*) \cap C(B_i^*)| = 0$.

If $a_i + b_i \geq \omega$, replicate $x_{2k}$ $\omega - b_i$ times and remove $x'_{2k}$. The new graph $H$ is
perfect and $\omega (H) = \omega$. Therefore $H$ has an $\omega$-coloring. Again $|C(A_i^*)|=a_i$ and
$|C(B_i^*)|=b_i$, and the $\omega - b_i$ colors that appear in the replicates of $x_{2k}$ belong to $C(A_i^*)$.
Since these colors cannot appear in $C(B_i^*)$, the number of common colors in $C(A_i^*)$ and $C(B_i^*)$ is $a_i +
b_i - \omega$.

(ii) $P_i$ has an even number of edges.

Assume w.l.o.g. that $a_i \leq b_i$. Let $P_i = x_1, \ldots, x_{2k+1}$. In $G_i^*$, replicate vertex $x_i$ $\omega
- a_i$ times for $i$ odd, $1 \leq i \leq 2k-1$, and replicate vertex $x_i$ $a_i$ times for $i$ even, $1 < i \leq
2k$. Finally, replicate $x_{2k+1}$ $\omega - b_i$ times. By Claim~1, the new graph $H$ is perfect and $\omega (H)
= \omega$. In an $\omega$-coloring of $H$, $|C(A_i^*)|=a_i$ and $|C(B_i^*)|=b_i$ and the colors that appear in
the replicates of $x_{2k}$ are precisely $C(A_i^*)$. But then these colors do not appear in the replicates of
$x_{2k+1}$ and consequently they must appear in $C(B_i)$. Thus $|C(A_i) \cap C(B_i)| = \min (a_i, b_i )$. This
proves Claim~3.
\\

\noindent{\bf Claim 4:} {\it $G^*$ has an $\omega
(G^*)$-coloring.}
\\

\noindent{\bf Proof of Claim~4:} Let $\omega = \omega (G^*)$.
Clearly, $\omega \geq a_1 + a_2$ and $\omega \geq b_1 + b_2$. To
prove the claim, we will combine $\omega$-colorings of $V_1^*$
and $V_2^*$.

If at least one of the sets $A_1^*, A_2^*, B_1^*, B_2^*$ is empty, one can easily construct the desired
$\omega$-coloring of $G^*$. So we assume now that these sets are nonempty. This implies that $\omega \geq \omega
(G_1^*)$ and $\omega \geq \omega (G_2^*)$. By Claim~3, there exist $\omega$-colorings of $V_i^*$ such that
$|C(A_i^*)| = a_i$ and $|C(B_i^*)| = b_i$. Thus, if $A_2^*$ and $B_2^*$ are in different connected components of
$G(V_2^*)$, an $\omega$-coloring of $V_1^*$ can be combined with $\omega$-colorings of the components of
$G(V_2^*)$ into an $\omega$-coloring of $G^*$. So we can assume that both $P_1$ and $P_2$ exist. Since $G_1$
contains no odd hole, every chordless path from $A_1$ to $B_1$ has the same parity as $P_1$. It follows from the
definition of 2-join decomposition that $P_1$ and $P_2$ have the same parity.

(i) $P_1$ and $P_2$ both have an odd number of edges.

Then by Claim~3 (i), there exists an $\omega$-coloring of $V_i^*$ with $|C(A_i^*) \cap C(B_i^*)| = \max (0, a_i +
b_i - \omega )$. In the coloring of $V_1^*$, label by 1 through $a_1$ the colors that occur in $A_1^*$ and by
$\omega$ through $\omega - b_1 +1$ the colors that occur in $B_1^*$. In the coloring of $V_2^*$, label by $\omega$
through $\omega - a_2 +1$ the colors that occur in $A_2^*$ and by 1 through $b_2$ the colors that occur in
$B_2^*$. If this is not an $\omega$-coloring of $G^*$, there must exist a common color in $A_1^*$ and $A_2^*$ or
in $B_1^*$ and $B_2^*$. But then either $a_1 \geq \omega - a_2 +1$ or $b_2 \geq \omega - b_1 +1$, a contradiction.

(ii) $P_1$ and $P_2$ both have an even number of edges.

Then by Claim~3 (ii), there exists an $\omega$-coloring of $V_i^*$ with $|C(A_i^*) \cap C(B_i^*)| = \min
(a_i,b_i)$. In the coloring of $V_1^*$, label by 1 through $a_1$ the colors that occur in $A_1^*$ and by 1 through
$b_1$ the colors that occur in $B_1^*$. In the coloring of $V_2^*$, label by $\omega$ through $\omega - a_2 +1$
the colors that occur in $A_2^*$ and by $\omega$ through $\omega - b_2 +1$ the colors that occur in $B_2^*$. If
this is not an $\omega$-coloring of $G$, there must exist a common color in $A_1^*$ and $A_2^*$ or in $B_1^*$ and
$B_2^*$. But then either $a_1 \geq \omega - a_2 +1$ or $b_1 \geq \omega - b_2 +1$, a contradiction.

\begin{CO} \label{2join}
If a minimally imperfect graph $G$ has a 2-join, then $G$ is an odd hole.
\end{CO}

\noindent {\bf Proof:}  Since $G$ is not perfect,
Theorem~\ref{2j} implies that block $G_1$ or $G_2$ is not
perfect, say $G_1$. Since $G_1$ is an induced subgraph of $G$ and
$G$ is minimally imperfect, it follows that $G=G_1$. Since $|V_2|
\geq 3$, $V_2$ induces a chordless path. Thus $G$ is a minimally
imperfect graph with a vertex of degree 2. This implies that $G$
is an odd hole \cite{P74}.

We end this section with another decomposition that preserves perfection. A graph $G$ has a {\em 6-join} if $V(G)$
can be partitioned into eight nonempty sets $X_1, X_2, X_3, X_4, Y_1, Y_2, Y_3, Y_4$ with the property that, for
any $x_i \in X_i$ ($i=1,2,3$) and $y_j \in Y_j$ ($j=1,2,3$), the graph induced by $x_1, y_1, x_2, y_2, x_3, y_3$
is a 6-hole and these kinds of edges are the only adjacencies between $X= X_1 \cup X_2 \cup X_3 \cup X_4$ and $Y =
Y_1 \cup Y_2 \cup Y_3 \cup Y_4$.

\begin{THTH} \label{6join}
{\em (Aossey and Vu\v{s}kovi\'c  \cite{av})} No minimally imperfect graph contains a 6-join.
\end{THTH}

If $G$ contains a 6-join, define {\em blocks} $G_X$ and $G_Y$ as follows. $G_X$ is the graph induced by $X\cup \{
y_1, y_2, y_3 \}$ where $y_j \in Y_j$ ($j=1,2,3$). Similarly $G_Y$ is the graph induced by $Y\cup \{ x_1, x_2, x_3
\}$ where $x_i \in X_i$ ($i=1,2,3$). It can be shown \cite{ao} that $G$ is perfect if and only if its blocks $G_X$
and $G_Y$ are perfect.

\section*{4. Skew Partition and Homogeneous Pair} \addsec

\vskip-5mm \hspace{5mm}

A graph has a {\it skew partition} if its vertices can be partitioned into four nonempty sets $A, B, C, D$ such
that there are all the possible edges between $A$ and $B$ and no edges from $C$ to $D$. It is easy to verify that
the odd holes and their complements do not have a skew partition. Chv\'atal \cite{vasek} conjectured that no
minimally imperfect graph has a skew partition.

\begin{THTH} \label{skew} {\bf (Skew Partition Theorem)}
{\em (Chudnovsky and Seymour  \cite{chsey2})} No minimally imperfect graph has a skew partition.
\end{THTH}

Chudnovsky and Seymour obtained this result as a consequence of their proof of the SPGC. In order to prove the
SPGC, they first proved the following weaker result.

\begin{THTH} \label{skewminimum}
{\em (Chudnovsky and Seymour  \cite{chsey})} A minimally imperfect Berge graph with smallest number of vertices
does not have a skew partition.
\end{THTH}

We do not give the proof of this difficult theorem here. Instead, we prove results due to Ho\`{a}ng \cite{Hoang96}
on two special skew partitions called $T$-cutset and $U$-cutset respectively.

Assume that $G$ is a minimally imperfect graph with skew partition $A, B, C, D$. Let $a = \omega (A)$, $b = \omega
(B)$, $\omega = \omega (G)$ and $\alpha = \alpha (G)$. The vertex sets $A \cup B \cup C$ and $A \cup B \cup D$
induce perfect graphs $G_1$ and $G_2$ respectively and both of these graphs contain an $\omega$-clique. Indeed,
each vertex of a minimally imperfect graph belongs to $\omega$ $\omega$-cliques \cite{P74} and, for $u \in C$,
these $\omega$-cliques are contained in $G_1$. For $u \in D$, they are contained in $G_2$.

\begin{LE} \label{l1}
 {\em (Ho\`{a}ng \cite{Hoang96})}
Let ${\cal C}_i$ be an $\omega$-coloring of $G_i$, for $i=1,2$. Then ${\cal C}_1$ and ${\cal C}_2$ cannot have the
same number of colors in $A$.
\end{LE}

\noindent {\bf Proof:} Suppose ${\cal C}_1$ and ${\cal C}_2$ have
the same number of colors in $A$ and assume w.l.o.g. that these
colors are $1,2, \ldots, k$. Let $K$ be the subgraph of $G$
induced by the vertices with colors $1, 2, \ldots, k$ and let $H
= G \setminus K$. Since every $\omega$-clique of $G$ is in $G_1$
or $G_2$, the largest clique in $K$ has size $k$ and the largest
clique in $H$ has size $\omega - k$. The graphs $H$ and $K$ are
perfect since they are proper subgraphs of $G$. Color $K$ with
$k$ colors and $H$ with $\omega - k$ colors. Now $G$ is colored
with $\omega$ colors, a contradiction to the assumption that $G$
is minimally imperfect.

\begin{LE} \label{l2}
No $\omega$-clique is contained in $A \cup B$.
\end{LE}

\noindent {\bf Proof:} Suppose that a $\omega$-clique were
contained in $A \cup B$. Then any $\omega$-coloring of $G_i$, for
$i=1,2$, would contain $a$ colors in $A$ and $b =\omega - a$
colors in $B$, contradicting Lemma~\ref{l1}.

\begin{LE} \label{l3}
Every $\alpha$-stable set intersects $A \cup B$.
\end{LE}

\noindent {\bf Proof:} By Lemma~\ref{l2} applied to the complement
graph, no $\alpha$-stable set is contained in $C \cup D$.

\begin{LE} \label{3.5}
If some $u \in A$ has no neighbor in $C$, then there exists an $\omega$-coloring of $G_1$ with $b$ colors in $B$.
\end{LE}

\noindent {\bf proof:} Let ${\cal C}_1$ be an $\omega$-coloring
of $G_1$ with minimum number $k$ of colors in $B$ and suppose
that this number is strictly greater than $b$. Consider the
subgraph $H$ of $G_1$ induced by the vertices colored with the
colors of ${\cal C}_1$ that appear in $B$. The graph $H \cup u$
can be colored with $k$ colors since it is perfect and has no
clique of size greater than $k$. Keeping the other colors of
${\cal C}_1$ in $G_1 \setminus (H \cup u)$, we get an
$\omega$-coloring of $G_1$ with fewer colors on $B$ than ${\cal
C}_1$, a contradiction.

\begin{LE} \label{l4}
If some $u \in A$ has no neighbor in $C$, then every vertex of $A$ has a neighbor in $D$ and every vertex of $B$
has a neighbor in $C$.
\end{LE}

\noindent {\bf Proof:} By Lemma~\ref{3.5}, there exists an
$\omega$-coloring of $G_1$ with $b$ colors in $B$. Thus, by
Lemma~\ref{l1}, there exists no $\omega$-coloring of $G_2$ with
$b$ colors in $B$. By Lemma~\ref{3.5}, this implies that every
vertex of $A$ has a neighbor in $D$.

Suppose that $v \in B$ has no neighbor in $C$. In the complement graph, $u$ and $v$ are adjacent to all the
vertices of $C$. By Lemma~\ref{l1}, $|A| \geq 2$ and $|B| \geq 2$.
 So $A'=A-u$,  $B'=B-v$, $C'= C$, $D'=D \cup
\{ u, v \}$ form a skew partition. But $u$ has no neighbor in $B$ and $v$ has no neighbor in $A$, contradicting
the first part of the lemma. So every  $v \in B$ has a neighbor in $C$.

A {\it $T$-cutset} is a skew partition with $u \in C$ and $v \in D$ such that every vertex of $A$ is adjacent to
both $u$ and $v$.

\begin{LE} \label{l5}  {\em (Ho\`{a}ng \cite{Hoang96})}
No minimally imperfect graph contains a $T$-cutset.
\end{LE}

\noindent {\bf Proof:}  In the complement, $u$ and $v$ contradict
Lemma~\ref{l4}.

A {\it $U$-cutset} is a skew partition with $u,v \in C$ such that every vertex of $A$ is adjacent to $u$ and every
vertex of $B$ is adjacent to $v$.

\begin{LE} \label{l6}  {\em (Ho\`{a}ng \cite{Hoang96})}
No minimally imperfect graph contains a $U$-cutset.
\end{LE}

\noindent {\bf Proof:}  In the complement, $u$ and $v$ contradict
Lemma~\ref{l4}.

We conclude this section with the notion of homogeneous pair introduced by Chv\'atal and Sbihi \cite{CS}. A graph
$G$ has a {\em homogeneous pair} if $V(G)$ can be partitioned into subsets $A_1$, $A_2$ and $B$, such that:

\begin{itemize}
\item  $|A_1| + |A_2| \geq 3$ and $|B|\geq 2$.

\item  If a node of $B$ is adjacent to a node of $A_1$ ($A_2$) then it is
adjacent to all the nodes of $A_1$ ($A_2$).
\end{itemize}

\begin{THTH} \label{Mjoin}
{\em (Chv\'atal and Sbihi \cite{CS})} No minimally imperfect graph contains a homogeneous pair.
\end{THTH}

\section*{5. Decomposition of Berge Graphs} \addsec

\vskip-5mm \hspace{5mm}

A graph is a {\em Berge graph} if it does not contain an odd hole or its complement. Clearly, all perfect graphs
are Berge graphs. The SPGC states that the converse is also true.

\begin{CN}  {\bf (Decomposition Conjecture)} \label{BGDC}
{\em (Conforti, Cornu\'ejols,
Robertson, Seymour, Thomas and Vu\v{s}kovi\'c (2001))} Every
Berge graph $G$ is basic or has a skew partition or a homogeneous
pair, or $G$ or $\bar G$ has a 2-join.
\end{CN}

This conjecture implies the SPGC. Indeed, suppose that the Decomposition Conjecture holds but not the SPGC. Then
there exists a minimally imperfect graph $G$ distinct from an odd hole or its complement. Choose $G$ with the
smallest number of vertices. $G$ is a Berge graph and it cannot have a skew partition by
Theorem~\ref{skewminimum}. $G$ cannot have an homogeneous pair by Theorem~\ref{Mjoin}. Neither $G$ nor $\bar G$
can have a 2-join by Corollary~\ref{2join}. So $G$ must be basic by the Decomposition Conjecture. Therefore $G$ is
perfect, a contradiction.

Note that there are other decompositions that cannot occur in minimally imperfect Berge graphs, such as 6-joins
(Theorem~\ref{6join}) or universal 2-amalgams \cite{CCGV} (universal 2-amalgams generalize both 2-joins and
homogeneous pairs). These decompositions could be added to the statement of Conjecture~\ref{BGDC} while still
implying the SPGC. However they do not appear to be needed. Paul Seymour commented that homogeneous pairs might
not be necessary either. In fact, we had initially formulated Conjecture~\ref{BGDC} without homogeneous pairs. I
added them to the statement to be on the safe side since they currently come up in the proof of the SPGC (see
below).

Several special cases of Conjecture~\ref{BGDC} are known. For example, it holds when $G$ is a Meyniel graph
(Burlet and Fonlupt \cite{BUR} in 1984), when $G$ is claw-free (Chvatal and Sbihi \cite{CS88} in 1988 and Maffray
and Reed \cite{MR} in 1999), diamond-free (Fonlupt and Zemirline \cite{FZ} in 1987), bull-free (Chv\'atal and
Sbihi \cite{CS} in 1987), or dart-free (Chv\'atal, Fonlupt, Sun and Zemirline \cite{CFSZ} in 2000). All these
results involve special types of skew partitions (such as star cutsets) and, in some cases, homogeneous pairs
\cite{CS}. A special case of 2-join called augmentation of a flat edge appears in  \cite{MR}. In 1999, Conforti
and Cornu\'ejols \cite{CC99} used more general 2-joins to prove Conjecture~\ref{BGDC} for WP-free Berge graphs, a
class of graphs that contains all bipartite graphs and all line graphs of bipartite graphs. This paper was the
precursor of a sequence of decomposition results involving 2-joins:

\begin{THTH} \label{sq}
{\em (Conforti, Cornu\'ejols and  Vu\v{s}kovi\'c \cite{square})} A square-free Berge graph is bipartite, the line
graph of a bipartite graph, or has a 2-join or a star cutset.
\end{THTH}

\begin{THTH}{\em (Chudnovsky, Robertson, Seymour and Thomas \cite{CRST})}
\label{3c} If $G$ is a Berge graph that contains the line graph of a bipartite subdivision of a 3-connected graph,
then $G$ has a skew partition, or $G$ or $\bar G$ has a 2-join or is the line graph of a bipartite graph.
\end{THTH}

Given two vertex disjoint triangles $a_1,a_2,a_3$ and $b_1,b_2,b_3$, a {\em stretcher} is a graph induced by three
chordless paths, $P^1=a_1, \ldots ,b_1$, $P^2=a_2, \ldots ,b_2$ and $P^3=a_3, \ldots ,b_3$, at least one of which
has length greater than one, such that $P^1, P^2, P^3$ have no common vertices and the only adjacencies between
the vertices of distinct paths are the edges of the two triangles. The next result is a real tour-de-force and a
key step in the proof of the SPGC.

\begin{THTH}{\em (Chudnovsky and Seymour \cite{chsey})}
\label{str} If $G$ is a Berge graph that contains a stretcher, then $G$ is the line graph of a bipartite graph or
$G$ has a skew partition or a homogeneous pair, or $G$ or $\bar G$ has a 2-join.
\end{THTH}

 A {\em wheel} $(H,v)$ consists of a hole $H$ together with a vertex $v$,
called the {\em center}, with at least three neighbors in $H$. If $v$ has $k$ neighbors in $H$, the wheel is
called a {\em $k$-wheel}. A {\em line wheel} is a 4-wheel $(H,v)$ that contains exactly two triangles and these
two triangles have only the center $v$ in common. A {\em twin wheel} is a 3-wheel containing exactly two
triangles. A {\em universal wheel} is a wheel $(H,v)$ where the center $v$ is adjacent to all the vertices of $H$.
A {\em triangle-free wheel} is a wheel containing no triangle. A {\em proper wheel} is a wheel that is not any of
the above four types. These concepts were first introduced in \cite{CC99}. The following theorem generalizes an
earlier result by Conforti, Cornu\'ejols and Zambelli \cite{CCZ} and Thomas \cite{thomas}.

\begin{THTH}{\em (Conforti, Cornu\'ejols and Zambelli \cite{CZ})}
\label{noproperw} If $G$ is a Berge graph that contains no proper wheels, stretchers or their complements, then
$G$ is basic or has a skew partition.
\end{THTH}

The last step in proving the SPGC is the following difficult theorem.

\begin{THTH}{\em (Chudnovsky and Seymour \cite{chsey2})}
\label{properw} If $G$ is a Berge graph that contains a proper wheel, but no stretchers or their complements, then
$G$ has a skew partition, or $G$ or $\bar G$ has a 2-join.
\end{THTH}

A monumental paper containing these results is forthcoming \cite{CRST2}. Independently, Conforti, Cornu\'{e}jols,
Vu\v{s}kovi\'{c} and Zambelli \cite{CCVZ} proved that the Decomposition Conjecture holds for Berge graphs
containing a large class of proper wheels but, as of May 2002, they could not prove it for all proper wheels.
Theorems~\ref{str}, \ref{noproperw} and \ref{properw} imply that Conjecture~\ref{BGDC} holds, and therefore the
SPGC is true.

\begin{CO}
{\bf (Strong Perfect Graph Theorem)} The only minimally imperfect graphs are the odd holes and their complements.
\end{CO}

Conforti, Cornu\'{e}jols and Vu\v{s}kovi\'{c}  \cite{doublestar} proved a weaker version of the Decomposition
Conjecture where ``skew partition'' is replaced by ``double star cutset''. A {\em double star} is a vertex set $S$
that contains two adjacent vertices $u, v$ and a subset of the vertices adjacent to $u$ or $v$. Clearly, if $G$
has a skew partition, then $G$ has a double star cutset: Take $S = A \cup B$,
 $u \in A$ and $v \in B$. Although the decomposition result in
\cite{doublestar} is weaker than Conjecture~\ref{BGDC} for Berge graphs, it holds for a larger class of graphs
than Berge graphs: By changing the decomposition from ``skew partition'' to ``double star cutset'', the result can
be obtained for all odd-hole-free graphs instead of just Berge graphs.

\begin{THTH} \label{odddecomp}
{\em (Conforti, Cornu\'{e}jols and Vu\v{s}kovi\'{c}  \cite{doublestar})} If $G$ is an odd-hole-free graph, then
$G$ is a bipartite graph or the line graph of a bipartite graph or the complement of the line graph of a bipartite
graph, or $G$ has a double star cutset or a 2-join.
\end{THTH}

One might try to use Theorem \ref{odddecomp}  to construct a polynomial time recognition algorithm  for
odd-hole-free graphs. Conforti, Cornu\'ejols, Kapoor and Vu\v{s}kovi\'c \cite{cckv3} obtained a polynomial time
recognition algorithm for the class of even-hole-free graphs. This algorithm is based on the decomposition of
even-hole-free graphs by 2-joins, double star and triple star cutsets obtained in \cite{cckv2}.

A useful tool for studying Berge graphs is due to Roussel and Rubio \cite{RR}. This lemma was proved independently
by Robertson, Seymour and Thomas \cite{RST}, who popularized it and named it {\it The Wonderful Lemma}. It is used
repeatedly in the proofs of Theorems~\ref{3c}-\ref{properw}.

\begin{LE} {\bf (The Wonderful Lemma)}
{\em (Roussel and Rubio \cite{RR})} \label{lemma:wonder} Let $G$ be a Berge graph and assume that $V(G)$ can be
partitioned into a set $S$ and an odd chordless path $P=u,u',\ldots,v',v$ of length at least $3$ such that $u$,
$v$ are both adjacent to all the vertices in $S$ and $\bar{G} (S)$ is connected. Then one of the following holds:
\begin{itemize}
\item[(i)] An odd number of edges of $P$ have both ends adjacent to
all the vertices in $S$.
\item[(ii)] $P$ has length 3 and $\bar{G} (S\cup \{u',v'\})$ contains an odd
chordless path between $u'$ and $v'$.
\item[(iii)] $P$ has length at least 5 and there exist two
nonadjacent vertices $x$, $x'$ in $S$ such that $(V(P)\sm\{u,v\})\cup\{x,x'\}$ induces a path.
\end{itemize}
\end{LE}

\label{lastpage}

\end {document}